\documentclass[12pt]{amsart}

\usepackage{graphicx}
\newcommand{\Z}{\Bbb Z}
\newcommand{\R}{\Bbb R}
\newcommand{\C}{\Bbb C}

\newtheorem{theorem}{Theorem}[section]
\newtheorem{proposition}[theorem]{Proposition}
\newtheorem{lemma}[theorem]{Lemma}

\newtheorem{corollary}[theorem]{Corollary}

\theoremstyle{definition}
\newtheorem{definition}[theorem]{Definition}
\newtheorem{example}[theorem]{Example}

\theoremstyle{remark}
\newtheorem{remark}[theorem]{Remark}

\numberwithin{equation}{section}

\newcommand{\aaa}{\mbox{$\alpha$}}
\newcommand{\bbb}{\mbox{$\beta$}}

\newcommand{\ggg}{\mbox{$\gamma$}}

\newcommand{\bdd}{\mbox{$\partial$}}

\newcommand{\T}{\widetilde}
\newcommand{\p}{\partial}

\begin{document}

\title[Floer homology of $(1,1)$-knots]{Knot Floer homology of $(1,1)$-knots}

\author[Goda, Matsuda and Morifuji]{Hiroshi Goda,
Hiroshi Matsuda and Takayuki Morifuji}

\address{
Department of Mathematics,
Tokyo University of Agriculture and Technology,
2-24-16, Naka, Koganei, Tokyo 184-8588, Japan}

\email{goda@cc.tuat.ac.jp, morifuji@cc.tuat.ac.jp}

\address{
Department of Mathematics, Hiroshima University, 
1-3-1 Kagamiyama, Higashi-Hiroshima city, Hiroshima 739-8526, Japan}

\email{matsuda@ms.u-tokyo.ac.jp}

\begin{abstract}
We present a combinatorial method
for a calculation of knot Floer homology
with $\Z$-coefficient of $(1, 1)$-knots,
and then demonstrate it for non-alternating
$(1,1)$-knots with ten crossings and the
pretzel knots of type $(-2,m,n)$.
Our calculations determine
the unknotting numbers and 4-genera of
the pretzel knots of this type.
\end{abstract}

\thanks{
The first and third authors are partially supported
by Grant-in-Aid for Scientific Research, (No. 15740031 and 14740036),
Ministry of Education, Science, Sports and Technology, Japan.
The second author is partially supported
by JSPS Research Fellowships for Young Scientists.}

\maketitle

\section{Introduction}

Knot Floer homology for knots in closed 3-manifolds
is defined by P. Ozsv{\' a}th and Z. Szab{\' o} in \cite{knot}.
In the definition,
they use Heegaard Floer homology \cite{os1}
for closed 3-manifolds.
It is known that 
an invariant of contact structures on closed 3-manifolds
\cite{contact}, and an invariant of closed 4-manifolds \cite{4-mfd} are
obtained from Heegaard Floer homology.
The estimates for the genus of knots \cite{knot},
the 4-genus and the unknotting number \cite{4-genus}
are also known to be obtained from 
knot Floer homology.

In general,
it is difficult to calculate Floer homology explicitly,
however,
Heegaard Floer homologies for lens spaces
and Seifert fibered spaces
are calculated in \cite[Proposition 3.1]{os2}
and \cite{os3} respectively.
Furthermore 
Heegaard Floer homologies for 3-manifolds
obtained by Dehn surgeries on
2-bridge knots in $S^3$ are
calculated by Rasmussen \cite{rasmussen},
using the method presented in
\cite[Proposition 3.2]{os2}.
Ozsv{\' a}th and Szab{\' o} presented in \cite{alternating}
a method for a calculation of knot Floer homology
for alternating knots in $S^3$, 
and they also calculated in \cite{genus}
knot Floer homology for knots with 
at most nine crossings, Kinoshita-Terasaka
knots and Conway knots.

In this paper,
we present a combinatorial method
for a calculation of knot Floer homology
with $\Z$-coefficient
which can be applied to all $(1, 1)$-knots, 
that is, all knots which admit $(1, 1)$-decompositions in $S^3$.
These knots form
one of wide and important classes in knot theory.
In fact,
it is well-known that torus knots and
2-bridge knots are a proper subset of $(1, 1)$-knots, 
and there exist $(1, 1)$-knots
which are non-alternating, 
or have arbitrary large crossing number.

Historically, 
Doll introduced in \cite{doll} 
the notion of $(g,b)$-decompositions
of knots and links in closed orientable 3-manifolds.
This is a generalization of $b$-bridge
decompositions of links in $S^3$. 
In this point of view, 
a $b$-bridge decomposition of a link in $S^3$ 
just corresponds to a $(0,b)$-decomposition. 
In this paper 
we focus on the case of $g=1$ and $b=1$, 
called $(1,1)$-knots.
A precise definition of $(1,1)$-decompositions of knots in $S^3$
is given in Section 2.
Recently 
$(1, 1)$-knots are extensively studied. 
See for example, \cite{berge}, \cite{c-m2}, \cite{choi-ko1},
\cite{goda-hayashi}, \cite{g-h-s},
\cite{hayashi1},
\cite{matsuda}, \cite{mulazzani},
\cite{munoz}, \cite{saito} and 
\cite{sch}.

This paper is organized as follows.
In the next section,
we briefly recall the notion of $(1,1)$-knots in $S^3$
and give their genus $2$ Heegaard splitting explicitly.
In Section 3,
we illustrate how to calculate 
the knot Floer homology
$\widehat{HFK}(S^3,K,i)$ through a sample calculation
for $K=10_{161}$.
In particular,
we explicitly determine the sign 
of each term appeared in the boundary operator 
of the complex $CFK^\infty(S^3,K)$. 
In Section 4,
we give a list of $\widehat{HFK}$
for non-alternating $(1,1)$-knots with ten crossings. 
Here, we present two knots whose knot Floer homologies 
are completely same (Example 4.4). 
In the final section,
we explicitly calculate $\widehat{HFK}$
for the pretzel knot of type $(-2,m,n)$,
where $m$ and $n$ are positive odd integers.
The result can be described by using the genus of them.

A part of this work was carried out while the first author
was visiting at Max-Planck-Institut f\"ur mathematik at Bonn.
He would like to express his sincere thanks for their hospitality.


\section{A genus 2 Heegaard splitting
for the complements of $(1,1)$-knots}

\begin{definition}
Let $V$ be a solid torus.
A properly embedded arc $t$ in $V$ is {\it trivial\/}
if there exists a disc $C$ embedded in $V$ satisfying
the following two properties:

\begin{enumerate}

\item $C \cap t = t$ is a subarc of $\partial C$;

\item $C \cap \partial V = c\ell (\partial C - t; \partial C)$
is an arc connecting the two points of $\partial t$ on $\partial V$.

\end{enumerate}

\end{definition}

\begin{definition}
A knot $K$ in $S^3$ is a {\it $(1,1)$-knot\/}
if there exists a genus 1 Heegaard splitting
$V_\alpha \cup V_\beta$ of $S^3$
satisfying the following two properties:

\begin{enumerate}

\item 
$K$ intersects the torus
$\partial V_\alpha = \partial V_\beta$
transversely in two points;

\item
the pair $(V_\alpha, K \cap V_\alpha)$
(resp. $(V_\beta, K \cap V_\beta)$) is
a pair of a solid torus and one trivial arc
properly embedded in $V_\alpha$ (resp. $V_\beta$).

\end{enumerate}

The decomposition $(V_\alpha, K \cap V_\alpha) \cup
(V_\beta, K \cap V_\beta)$ is called a {\it $(1,1)$-decomposition\/}
of $(S^3, K)$.
Let $\ggg$ be an embedded arc in $V_{\aaa}$ 
such that 
$\bdd\ggg\subset\text{int}(K\cap V_\alpha)$. 
Then $K\cap V_\alpha$ is cut into three arcs. 
The boundaries of one of them coincide with those of $\ggg$, 
call $\delta$. If the boundary of the regular neighborhood 
of $\ggg\cup\delta$ in $V_{\aaa}$ is isotopic to $\bdd V_{\aaa}$, 
then we call $\ggg$ {\it $(1,1)$-tunnel\/} of $K$.
Conversely, if a knot $K$ has a $(1,1)$-tunnel,  
it induces a $(1,1)$-decomposition of $K$. 
\end{definition}

Let $C_\alpha$ (resp. $C_\beta$) be a disc
which realizes the triviality of the arc $K \cap V_\alpha$
(resp. $K \cap V_\beta$) in $V_\alpha$ (resp. $V_\beta$).
We remark here that
if $C_\alpha \cap C_\beta = \emptyset$
on $\partial V_\alpha = \partial V_\beta$,
then $K$ is a torus knot in $S^3$.
Let $D_\alpha^2$ (resp. $D_\beta^2$) be a meridian disc of
the solid torus $V_\alpha$ (resp. $V_\beta$)
which is disjoint from the disc $C_\alpha$ (resp. $C_\beta$).

\begin{proposition} \label{lem: 11}
Suppose that $K$ is a $(1,1)$-knot in $S^3$.
Then we can construct a genus $2$ Heegaard splitting
$W_\alpha \cup W_\beta$ of $S^3$, and
a meridian disc system $D_\alpha^1$, $D_\alpha^2$
$($resp. $D_\beta^1$, $D_\beta^2)$ of $W_\alpha$ 
$($resp. $W_\beta)$
satisfying the following three properties:

\begin{enumerate}

\item
$K$ is contained in the interior of $W_\beta$;

\item
$K$ intersects $D_\beta^1$ transversely
in exactly one point, and $K$ is disjoint from $D_\beta^2$;

\item
$\partial D_\beta^1$ intersects $\partial D_\alpha^1$
transversely in exactly one point
on $\partial W_\alpha = \partial W_\beta$, and
$\partial D_\beta^1$ is disjoint from $\partial D_\alpha^2$
on $\partial W_\alpha = \partial W_\beta$.
\end{enumerate}
\end{proposition}

\begin{remark}
The properties stated in the above lemma
satisfies the supposition of Proposition 6.1 in \cite{knot}. 
Thus 
we can directly apply a method discussed there 
for our situation. 
\end{remark}

\begin{proof}
Suppose that $K$ is a $(1,1)$-knot in $S^3$, and that
$(V_\alpha, t_\alpha) \cup (V_\beta, t_\beta)$ is
a $(1,1)$-decomposition of the pair $(S^3, K)$.
Drilling $V_\alpha$ along the arc $t_\alpha$,
we obtain a genus 2 Heegaard splitting $W_\alpha \cup W_\beta$
of $S^3$, where $W_\alpha$ is obtained from $V_\alpha$
by removing $N(t_\alpha; V_\alpha)$, and
$W_\beta$ is obtained from $V_\beta$ by attaching
$N(t_\alpha; V_\alpha)$ as a 1-handle.
We remark that $K$ is contained in the interior of $W_\beta$.
Let $D_\beta^1$ be a cocore disc of the 1-handle
$N(t_\alpha; V_\alpha)$.
Note that $D_\beta^1$ is a meridian disc of $W_\beta$ such that
$K$ intersects $D_\beta^1$ transversely
in exactly one point.
Let $D_\alpha^1$ be a disc in $W_\alpha$ which is
the restriction of the disc $C_\alpha$ in $V_\alpha$.
This disc $D_\alpha^1$ is a meridian disc of $W_\alpha$.
The above construction shows that
these meridian discs $D_\alpha^1$ and $D_\beta^1$
of $W_\alpha$ and $W_\beta$, respectively,
satisfy the property that
$\partial D_\alpha^1$ intersects $\partial D_\beta^1$
transversely in exactly one point
on $\partial W_\alpha = \partial W_\beta$.
We denote by $D_\alpha^2$ (resp. $D_\beta^2$) the image
in $W_\alpha$ (resp. $W_\beta$)
of the meridian disc $D_\alpha^2$ (resp. $D_\beta^2$)
of $V_\alpha$ (resp. $V_\beta$).
Since $D_\beta^2$ is disjoint from the disc $C_\beta$ in $V_\beta$,
$K$ is disjoint from  $D_\beta^2$.
Since $D_\alpha^2$ is disjoint from the disc $C_\alpha$
in $V_\alpha$,
$\partial D_\beta^1$ is disjoint from $\partial D_\alpha^2$
on $\partial W_\alpha = \partial W_\beta$.
Note that $\{ D_\alpha^1, D_\alpha^2 \}$
(resp. $\{ D_\beta^1, D_\beta^2 \})$ is
a complete meridian disc system of
the genus 2 handlebody $W_\alpha$ (resp. $W_\beta$).
\end{proof}


\section{A sample calculation of $\widehat{HFK}$
for the knot $10_{161}$}

As for the definition of the knot Floer homology
$\widehat{HFK}(S^3,K,i)$ for a knot $K$ in $S^3$,
see the original paper \cite{knot}.
Throughout this paper,
we use the same notation in \cite{knot} 
and omit the explanation here. 

Let $K$ be the knot $10_{161}$ in Rolfsen's table \cite{rol} 
with its orientation indicated in the Figure $10_{161}(0)$.

\begin{figure}
\centering
\includegraphics[width=8cm,height=7cm]{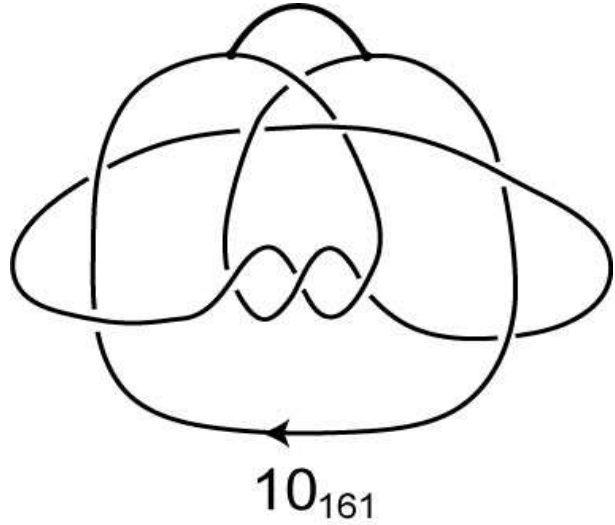}
\caption{$10_{161}(0)$}
\end{figure}

It is not difficult to check that 
this knot is a $(1,1)$-knot in $S^3$.
In fact, 
Figure $10_{161}(1)$ illustrates a $(1,1)$-decomposition
$(V_\alpha, t_\alpha) \cup (V_\beta, t_\beta)$ of $(S^3, K)$, and
the shaded disc in the figure illustrates a meridian disc
$D_\alpha^2$ of $V_\alpha$ which is disjoint from $C_\alpha$.
Further 
Figures $10_{161}(2) - (5)$ show that
the arc $t_\beta$ is trivial in $V_\beta$, and
they also show a construction of a meridian disc $D_\beta^2$
of $V_\beta$ which is disjoint from $C_\beta$.

\begin{figure}
\centering
\includegraphics[width=8cm,height=4.5cm]{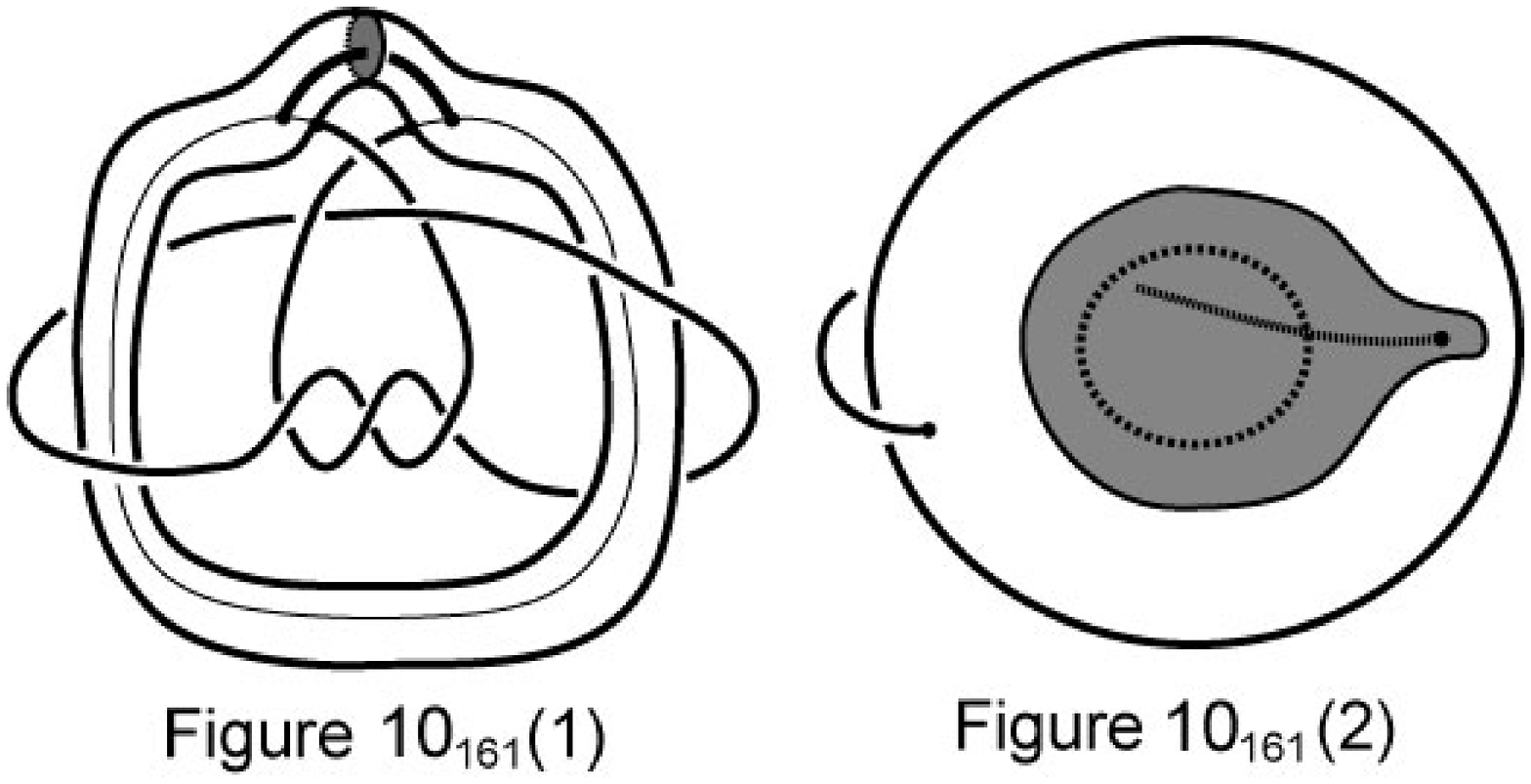}
\caption{}
\end{figure}

\begin{figure}
\centering
\includegraphics[width=8cm,height=4cm]{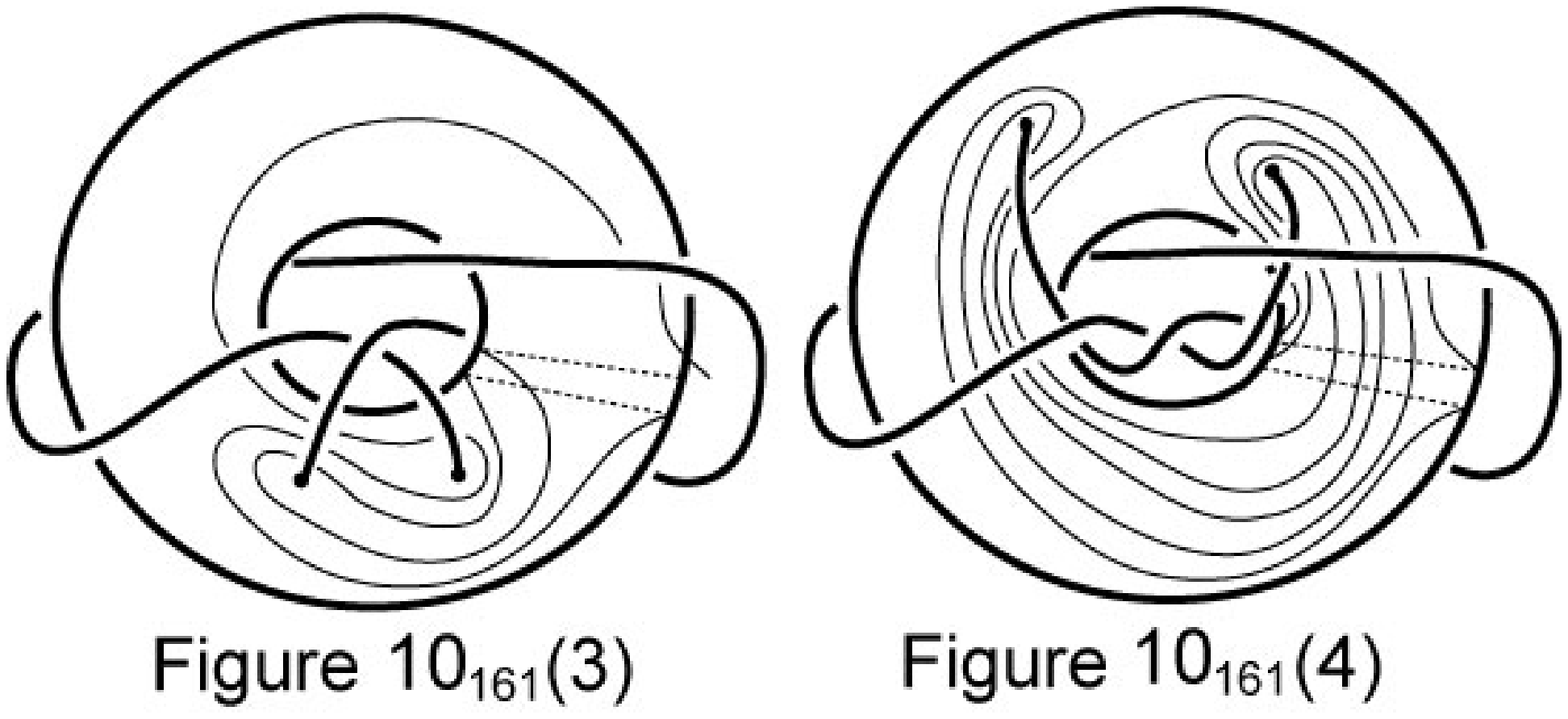}
\caption{}
\end{figure}

\begin{figure}
\centering
\includegraphics[width=7cm,height=7cm]{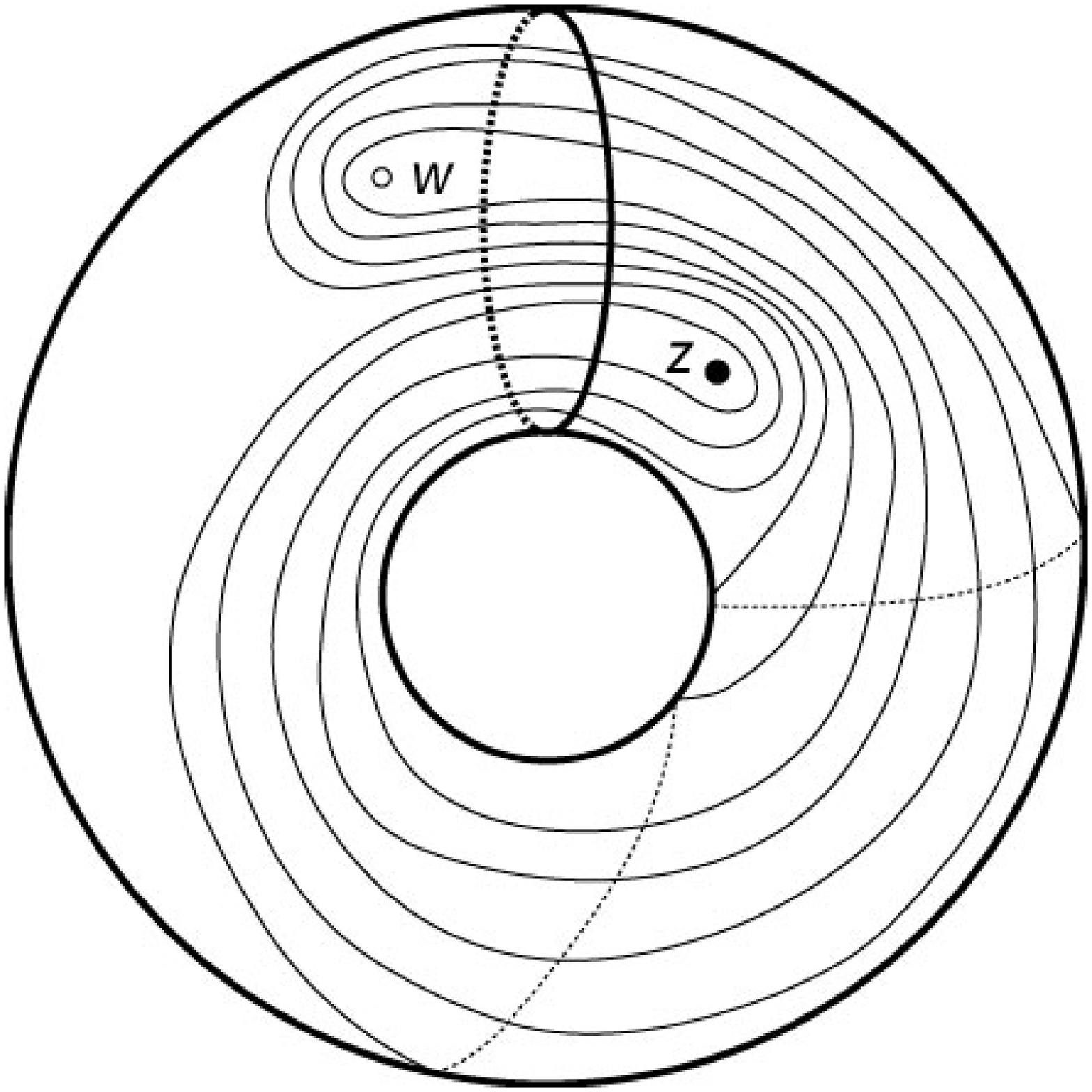}
\caption{$10_{161}(5)$}
\end{figure}

A genus $2$ Heegaard splitting of $S^3$
constructed in Proposition \ref{lem: 11}
can be destabilized along $D_\alpha^1$ and $D_\beta^1$.
In fact, 
Figure $10_{161}(5)$ illustrates a Heegaard diagram 
of genus 1 of $S^3$
which is obtained from that of genus 2
after the destabilization along $D_\alpha^1$ and $D_\beta^1$.
On the Heegaard surface of genus 1 
(which is a torus $\Sigma$), 
meridian curves $\partial D_\alpha^2$ and
$\partial D_\beta^2$, and two reference points $w$ and $z$
are illustrated in Figure $10_{161}(5)$.
Moreover, 
Figure $10_{161}(6)$ illustrates the corresponding diagram on the annulus
which is obtained from the torus $\Sigma$ by cutting along
$\partial D_\alpha^2$.

\begin{figure}
\centering
\includegraphics[width=7cm,height=5cm]{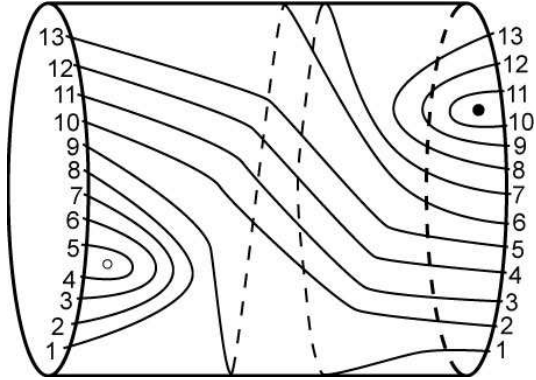}
\caption{$10_{161}(6)$}
\end{figure}

Finally, 
Figure $10_{161}(7)$ illustrates the corresponding diagram
of $\widetilde{\alpha}$ and $\widetilde{\beta}$
on the universal cover $\C$ of the torus $\Sigma$,
where $\widetilde{\alpha}$ (resp. $\widetilde{\beta}$)
is a lift of $\partial D_\alpha^2$ (resp. $\partial D_\beta^2$)
in $\C$.
Lightly-shaded (resp. darkly-shaded) circles represent
lifts of the reference point $w$ (resp. $z$) in $\C$.
The point labeled $\widetilde{x}$ on $\C$, 
which is an intersection point of $\widetilde{\alpha}$
and $\widetilde{\beta}$, 
is a lift of the point $x$ on the torus $\Sigma$.

\begin{figure}
\centering
\includegraphics[width=8.5cm,height=7.7cm]{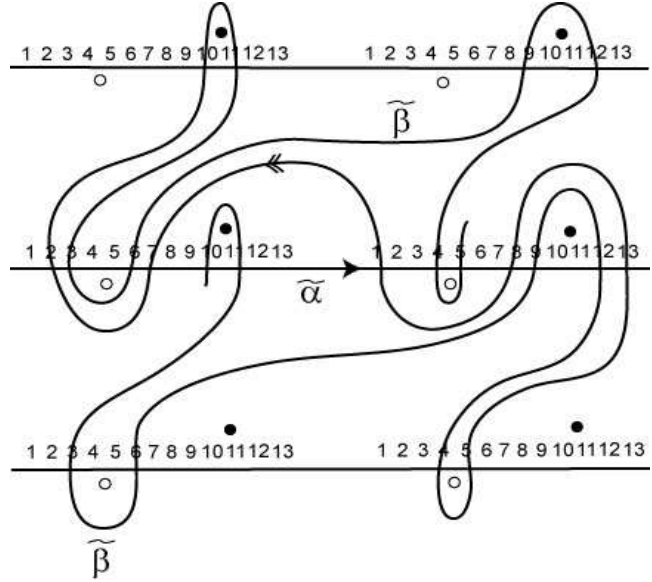}
\caption{$10_{161}(7)$}
\end{figure}

Proposition 6.4 in \cite{knot} shows that
if there exists a holomorphic disc
in $\pi_2 (\widetilde{x}, \widetilde{y})$,
then the absolute value of the coefficient of $\widetilde{y}$
in the expression of $\partial [\widetilde{x}; i,j]$ is 1.
In the following,
we try to 
explain how to assign $+1$ or $-1$
to the coefficient of $\widetilde{y}$
in the expression of $\partial [\widetilde{x}; i, j]$.

Let $u(x, y)$ denote a representative in $\C$ of
a holomorphic disc in $\pi_2 (\widetilde{x}, \widetilde{y})$,
where $\{ \widetilde{x}, \widetilde{y} \} \in
\{ \widetilde{x}_1, \widetilde{x}_2, \cdots, \widetilde{x}_{13} \}$.
Note that every holomorphic disc
in $\pi_2 (\widetilde{x}, \widetilde{y})$
can be represented as an embedded disc
in the universal cover $\C$ of $\Sigma$.
Then 
Figure $10_{161}(7)$ illustrates orientations of
the curves $\widetilde{\alpha}$ and $\widetilde{\beta}$
which are given as totally real submanifolds of $\C$.
Since the boundary of the disc $u(x, y)$ consists of
a subarc of $\widetilde{\alpha}$ and
a subarc of $\widetilde{\beta}$,
this orientation of $\widetilde{\alpha}$ induces
an orientation of the boundary of the disc $u(x, y)$, and
an orientation of the disc $u(x, y)$.

Now hereafter, 
we also use the same notation as in $\S$3.6 of \cite{os1},
but we use a connection over the universal cover $\C$
of ${\rm Sym}^1 (\Sigma) = \Sigma$
which is trivial along $\widetilde{\alpha}$. 
Next 
we choose a contraction of the disc $\mathbb{D}$
to preserve the subarc of $\partial \mathbb{D}$
mapping into $\widetilde{\alpha}$.
The linear transformation from
${\C} \cong {T_{u(0+i)}} \widetilde{\beta} \otimes_{\R} \C$
to $T_{u(0)} \widetilde{\beta} \otimes_{\R} \C \cong \C$
corresponds to a complex number $A(1) \in GL(1; \C) = \C$.
Then $A(1)$ induces an orientation of the disc $u(x, y)$.
We assign a sign
$\varepsilon (x, y) = +1$ (resp. $\varepsilon (x, y) = -1$)
to the disc $u(x, y)$
if the above two orientations of the disc $u(x, y)$
agree (resp. disagree).
These assignments $\varepsilon (x, y)$ give
an orientation of 1-dimensional subspaces of the moduli space.

In order to choose a coherent orientation for all moduli spaces,
we follow the arguments in $\S$21 of \cite{fooo}.
In this paragraph,
we use the same notation as in $\S$21 of \cite{fooo}.
The above choice of an orientation for 1-dimensional subspaces
of the moduli space gives an orientation for
${\rm Index} \, \bar{\partial}_{(u; \Lambda(x), \Lambda(y))}$.
We choose the same orientation on
${\rm Index} \, \bar{\partial}_{(-, \Lambda(x))}$
for each $\widetilde{x} \in
\widetilde{\alpha} \cap \widetilde{\beta}$, and
we choose the same orientation on
${\rm Index} \, \bar{\partial}_{(+, \Lambda(y))}$
for each $\widetilde{y} \in
\widetilde{\alpha} \cap \widetilde{\beta}$.
Then 
we see from Remark 21.12 (2), Proposition 23.2 and Lemma 23.4
in \cite{fooo} that
these choices of orientations give a coherent orientation for
all moduli spaces.
Hence we obtain the following proposition.

\begin{proposition}
Let $K$ be a $(1,1)$-knot in $S^3$.
Suppose that the corresponding Heegaard diagram on $\C$
is constructed as above.
Then the coefficient of the generator corresponding to the point
$\widetilde{y}$ in the expression of $\partial [\widetilde{x}; i, j]$
is equal to $\varepsilon (x, y)$.
\end{proposition}

For example, 
we assign $+1$ (resp. $-1$) to the coefficient
of the generator corresponding to the point $\widetilde{x}_1$
in the expression of $\partial [\widetilde{x}_7 ; i, j]$
(resp. $\partial [\widetilde{x}_8 ; i, j]$). 
From Figure $10_{161}(7)$, 
we can verify that 
the complex $CFK^\infty (S^3, K)$ is a $\Z$-module
generated by
\begin{align*}
& [\widetilde{x}_1; i, i],\ [\widetilde{x}_2; i-1, i+1],\
[\widetilde{x}_3; i-1, i],\ [\widetilde{x}_4; i, i-2],\
[\widetilde{x}_5; i+1, i-2],\\
& [\widetilde{x}_6; i,i],\
[\widetilde{x}_7; i, i+1],\
[\widetilde{x}_8; i+1, i],\ [\widetilde{x}_9; i, i],\
[\widetilde{x}_{10}; i-2, i+1],\\
& [\widetilde{x}_{11}; i-2, i],\
[\widetilde{x}_{12}; i, i-1]\
\mathrm{and}\ [\widetilde{x}_{13}; i+1, i-1]
\end{align*}
for $i \in \Z$, 
and the boundary operator on $CFK^\infty (S^3, K)$
is explicitly given by

\begin{align*}
\partial [\widetilde{x}_1; i, i] 
& = 
[\widetilde{x}_4; i, i-2] - [\widetilde{x}_{11}; i-2, i],\\
\partial [\widetilde{x}_2; i-1, i+1]
& =
[\widetilde{x}_3; i-1, i],\\
\partial [\widetilde{x}_3; i-1, i] 
& = 0,\\
\partial [\widetilde{x}_4; i, i-2] 
& = 0,\\
\partial [\widetilde{x}_5; i+1, i-2] 
&=  
- [\widetilde{x}_4; i, i-2],\\
\partial [\widetilde{x}_6; i, i] 
& = 
[\widetilde{x}_4; i, i-2] - [\widetilde{x}_3; i-1, i],\\
\partial [\widetilde{x}_7; i, i+1] 
& = 
[\widetilde{x}_1; i, i] - [\widetilde{x}_2; i-1, i+1]
- [\widetilde{x}_6; i, i] + [\widetilde{x}_{10}; i-2, i+1],\\
\partial [\widetilde{x}_8; i+1, i] 
& = 
[\widetilde{x}_9; i, i] + [\widetilde{x}_{13}; i+1, i-1]
- [\widetilde{x}_1; i, i] - [\widetilde{x}_5; i+1, i-2],\\
\partial [\widetilde{x}_9; i, i] 
& = 
[\widetilde{x}_{12}; i, i-1] - [\widetilde{x}_{11}; i-2, i],\\
\partial [\widetilde{x}_{10}; i-2, i+1] 
& = 
[\widetilde{x}_{11}; i-2, i],\\
\partial [\widetilde{x}_{11}; i-2, i] 
& = 0,\\
\partial [\widetilde{x}_{12}; i, i-1] 
& = 0,\\
\partial [\widetilde{x}_{13}; i+1, i-1] 
& = - [\widetilde{x}_{12}; i, i-1].
\end{align*}

The homology of the complex $CFK^{0, *} (S^3, K)$ is
${\Z} \cong \widehat{HF} (S^3)$,
which is generated by the cycle $[\widetilde{x}_5; 0, -3]$.
Its absolute grading is defined to be 0 (see \cite{4-mfd}). 
We then obtain the knot Floer homology of the knot
$10_{161}$ in $S^3$ as follows:

\begin{center}
$  {\widehat {HFK}} (S^3, 10_{161}, i) \cong
  \begin{cases}
    {{\Z}_{(6)}} & \text{$i = 3$} \\
    {{\Z}_{(4)}} \oplus {{\Z}_{(5)}} & \text{$i = 2$} \\
    {{\Z}_{(3)}^2} & \text{$i = 1$} \\
    {{\Z}_{(2)}^3} & \text{$i = 0$} \\
    {{\Z}_{(1)}^2} & \text{$i = -1$} \\
    {{\Z}_{(0)}} \oplus {{\Z}_{(1)}} & \text{$i = -2$} \\
    {{\Z}_{(0)}} & \text{$i = -3$} \\
    0 & \text{otherwise.}
  \end{cases} $
\end{center}

\begin{corollary}
The genus $g(K)$, the $4$-genus $g_4(K)$ and
the unknotting number $u(K)$ of the knot $K = 10_{161}$ are $3$.
\end{corollary}

\begin{proof}
The adjunction inequality, Theorem 5.1 in \cite{knot},
shows that $g(K) \geq 3$, and
it is easy to see that $K$ bounds a Seifert surface of genus three.
It follows that $g(K) = 3$.

We next calculate the invariant $\tau(K)$, defined in \cite{4-genus}.
The subcomplex ${\mathcal{F}} (K, -3)$ is generated by
the generator corresponding to the point $\widetilde{x}_5$, and
the homology $H_* ({\mathcal{F}} (K, -3))$ is generated by
the cycle corresponding to the point $\widetilde{x}_5$,
whose absolute grading is 0.
It follows that $\tau(K) = -3$.
Corollary 1.3 in \cite{4-genus} shows $g_4(K) \geq | \tau(K) | = 3$.
Since $g(K) = 3$, 
we have $g_4(K) \leq 3$ and 
hence $g_4(K) = 3$.
It is easy to see that crossing changes at three crossings
change $K$ to the trivial knot.
See, for example, \cite{4-genus} or \cite{tanaka}.
It is well-known that $u(K) \geq g_4(K)$, so we have $u(K) = 3$.
\end{proof}

\begin{remark}
The 4-genus and the unknotting number of the knot $10_{161}$
were determined independently by Tanaka \cite{tanaka}, 
Kawamura \cite{kawamura} and Ozsv{\' a}th-Szab{\' o} \cite{4-genus}.
\end{remark}


\section{Tunnel number 1 knots with 10 crossings}

Tunnel number 1 knots up to ten crossings are
determined in \cite{msy}.
In this section,
we present a list of the knot Floer homology
of non-alternating $(1,1)$-knots 
with ten crossings.

Knot Floer homology of alternating knots
is known to be described by using their signatures \cite{alternating}.
That of prime knots up to 
nine crossings are treated
in \cite{genus}.
We illustrate here
tunnel number 1 knots
up to ten crossings with $(1,1)$-tunnels 
except for these knots
at the end of this paper.
Note that the set of tunnel number 1 knots
include that of $(1,1)$-knots, and that
we can actually verify 
non-alternating tunnel number 1 knots
with ten crossings have $(1,1)$-decompositions. 
We use here the notation
in Rolfsen's book \cite{rol}.

\begin{lemma}
Non-alternating $(1,1)$-knots with ten crossings
are given in Figure $17$.
Moreover
$(1,1)$-tunnel of them are indicated by thick arcs
in each figure.
\end{lemma}

\begin{theorem}
Let $K=10_n$ be a non-alternating $(1,1)$-knot with ten crossings.
Then the knot Floer homology of $K$ 
and the invariant $\tau(K)$ are
given in Table $1$. 

\begin{table}[h]
\caption{$\widehat{HFK}(S^3,10_n,i)$}
\begin{center}
\begin{tabular}{|c||c||c|c|c|c|c|c|}
\noalign{\hrule height0.8pt}
$n$ & $\tau(K)$ & $i=0$ & $i=1$ & $i=2$ & $i=3$ & $i=4$ & $i\geq 5$ \\
\hline
$124$ &$4$& $\Z_{(-3)}$ & $\Z_{(-2)}$ & $0$ & $\Z_{(-1)}$ & $\Z_{(0)}$ & $0$ \\
\hline
$125$ &$1$&$\Z_{(-1)}$&$\Z_{(0)}^2$ & $\Z_{(1)}^2$ & $\Z_{(2)}$ & $0$ & $0$ \\
\hline
$126$ &$-1$&$\Z_{(1)}^5$&$\Z_{(2)}^4$ & $\Z_{(3)}^2$ & $\Z_{(4)}$ & $0$ & $0$ \\
\hline
$127$ &$-2$& $\Z_{(2)}^7$ & $\Z_{(3)}^6$ & $\Z_{(4)}^4$ & $\Z_{(5)}$ & $0$ & $0$ \\
\hline
$128$ & $3$&$\Z_{(-2)}$ & $\Z_{(-2)}$ & $\Z_{(-1)}^3$ & $\Z_{(0)}^2$ & $0$ & $0$ \\
\hline
$129$ &$0$& $\Z_{(0)}^9$ & $\Z_{(1)}^6$ & $\Z_{(2)}^2$ & $0$ & $0$ & $0$ \\
\hline
$130$ &$0$& $\Z_{(0)}^5$ & $\Z_{(1)}^4$ & $\Z_{(2)}^2$ & $0$ & $0$ & $0$ \\
\hline
$131$ &$-1$ & $\Z_{(1)}^{11}$ & $\Z_{(2)}^8$ & $\Z_{(3)}^2$ & $0$ & $0$ & $0$ \\
\hline
$132$ &$-1$& $\Z_{(0)}^2\oplus\Z_{(1)}$ & $\Z_{(1)}^2\oplus\Z_{(2)}$ &
        $\Z_{(2)}$ & $0$ & $0$ & $0$ \\
\hline
$133$ &$-1$& $\Z_{(1)}^7$ & $\Z_{(2)}^5$ & $\Z_{(3)}$ & $0$ & $0$ & $0$ \\
\hline
$134$ &$3$& $\Z_{(-3)}^3$ & $\Z_{(-2)}^4$ & $\Z_{(-1)}^4$ &
$\Z_{(0)}^2$ & $0$ & $0$ \\
\hline
$135$ & $0$&$\Z_{(0)}^{13}$ & $\Z_{(1)}^9$ & $\Z_{(2)}^3$ & $0$ & $0$ & $0$ \\
\hline
$136$ & $0$&$\Z_{(-1)}^6\oplus\Z_{(0)}$ & $\Z_{(0)}^4$ & $\Z_{(1)}$ &
$0$ & $0$ & $0$ \\
\hline
$137$ & $0$&$\Z_{(0)}^{11}$ & $\Z_{(1)}^6$ & $\Z_{(2)}$ & $0$ & $0$ & $0$ \\
\hline
$138$ &$1$& $\Z_{(-1)}^7$ & $\Z_{(0)}^8$ & $\Z_{(1)}^5$ & $\Z_{(2)}$ & $0$ & $0$ \\
\hline
$139$ & $4$&$\Z_{(-3)}^3$ & $\Z_{(-2)}^2$ & $0$ & $\Z_{(-1)}$ & $\Z_{(0)}$ & $0$ \\
\hline
$145$ & $-2$&$\Z_{(1)}^4\oplus\Z_{(2)}$ & $\Z_{(2)}^2\oplus\Z_{(3)}$ & $\Z_{(4)}$
      & $0$ & $0$ & $0$ \\
\hline
$161=162$ & $-3$&$\Z_{(2)}^3$ & $\Z_{(3)}^2$ & $\Z_{(4)}\oplus\Z_{(5)}$
          & $\Z_{(6)}$ & $0$ & $0$ \\
\noalign{\hrule height0.8pt}
\end{tabular}
\end{center}
\end{table}

\end{theorem}

\begin{remark}
As was stated in \cite{knot} Section 3,
it holds that
$$
\widehat{HFK}_d(S^3,K,i) \cong \widehat{HFK}_{d-2i}(S^3,K,-i).
$$
Thereby
we have given a list of $\widehat{HFK}(S^3,K,i)$
only for $i\geq 0$.
\end{remark}

\begin{proof}
We can directly apply the method discussed in Section 3 
to these $(1,1)$-knots. 
However 
the calculation is straightforward and lengthy, 
so that we omit here. 
\end{proof}

\begin{example}
Let us consider the knot $K_1=8_8$, 
which is an alternating $(1,1)$-knot. 
It has the Alexander polynomial 
$9-6(t+t^{-1})+2(t^2+t^{-2})$ and 
the signature zero. 
Thus 
we have $\widehat{HFK}(S^3,K_1,i)=\Z_{(i)}^{|a_i|}$, 
where $a_i$ denotes the coefficient of 
$i$-th term of the Alexander polynomial 
(see \cite{alternating}). 
Further, $\tau(K_1)=0$ from Theorem 1.4 in 
\cite{4-genus}.
We then see from Theorem 4.2 that 
the knot Floer homology and $\tau$ 
of $K_2=10_{129}$ 
are same as those of $K_1$. 
Further, 
it is known that these two knots have the 
same Alexander polynomial and Jone polynomial, 
but are not mutant each other \cite{L-M}. 
M. Teragaito informed us this example.
\end{example}





\section{Pretzel knots}

Tunnel number 1 pretzel knots are determined in
\cite{klimenko} and \cite{msy}. 
More precisely, 
the type is $(\pm 2,m,n)$, 
where $m$ and $n$ are odd integers. 
It is not difficult to check 
that all of them admit
$(1,1)$-decompositions.
In this section, we show the following theorem.
See \cite{kawauchi} 
as for the classification of
pretzel knots. 
In particular, 
we remark here that 
$P(\pm 2,m,n)$ is equivalent to $P(\pm 2,n,m)$.

\begin{theorem}
Let $K$ be the pretzel knot $P(-2,m,n)$,
where
$m$ and $n$ are odd integers
satisfying $m\geq n\geq 3$.
Then we have
\begin{center}
$
{\widehat{HFK}}(S^3,K,i) \ \cong
 \ \begin{cases}
    0 & \text{$(i>g)$} \\
    {\Z}_{(g+i)} & \text{$(i =g,\ g- 1)$} \\
    0 & \text{$(i = g-2)$} \\
    {\Z}^{g-2-i}_{(g-1+i)} & \text{$(g-n\leq i<g-2)$} \\
    {\Z}^{n-2}_{(g-1+i)} & \text{$(0\leq i < g-n)$}
  \end{cases}
$
\end{center}
where 
$\displaystyle{g=\frac{m+n}{2}}$ 
denotes the genus of the pretzel knot $K$.
Moreover, 
the invariant $\tau(K)$ is $-g$. 
\end{theorem}

\begin{remark}
Knot Floer homology of the
remaining tunnel number 1 pretzel knots
can be calculated by the same method
as in stated below.
\end{remark}

\begin{corollary}
The $4$-genus $g_4(K)$ and 
the unknotting number $u(K)$ of 
the pretzel knot $K=P(-2,m,n)$ are 
equal to the genus 
$\displaystyle{g=\frac{m+n}{2}}$. 
\end{corollary}

\begin{proof}
It is clear that 
the invariant $\tau(K)$ is equal to $-g$. 
Thereby 
the claim follows from the same argument 
in the proof of Corollary 3.2.
\end{proof}

At first, 
we present how to construct 
a genus $1$ Heegaard diagram
for pretzel knots with $(1,1)$-decompositions.

\begin{figure}
\centering
\includegraphics[width=8.3cm,height=4.7cm]{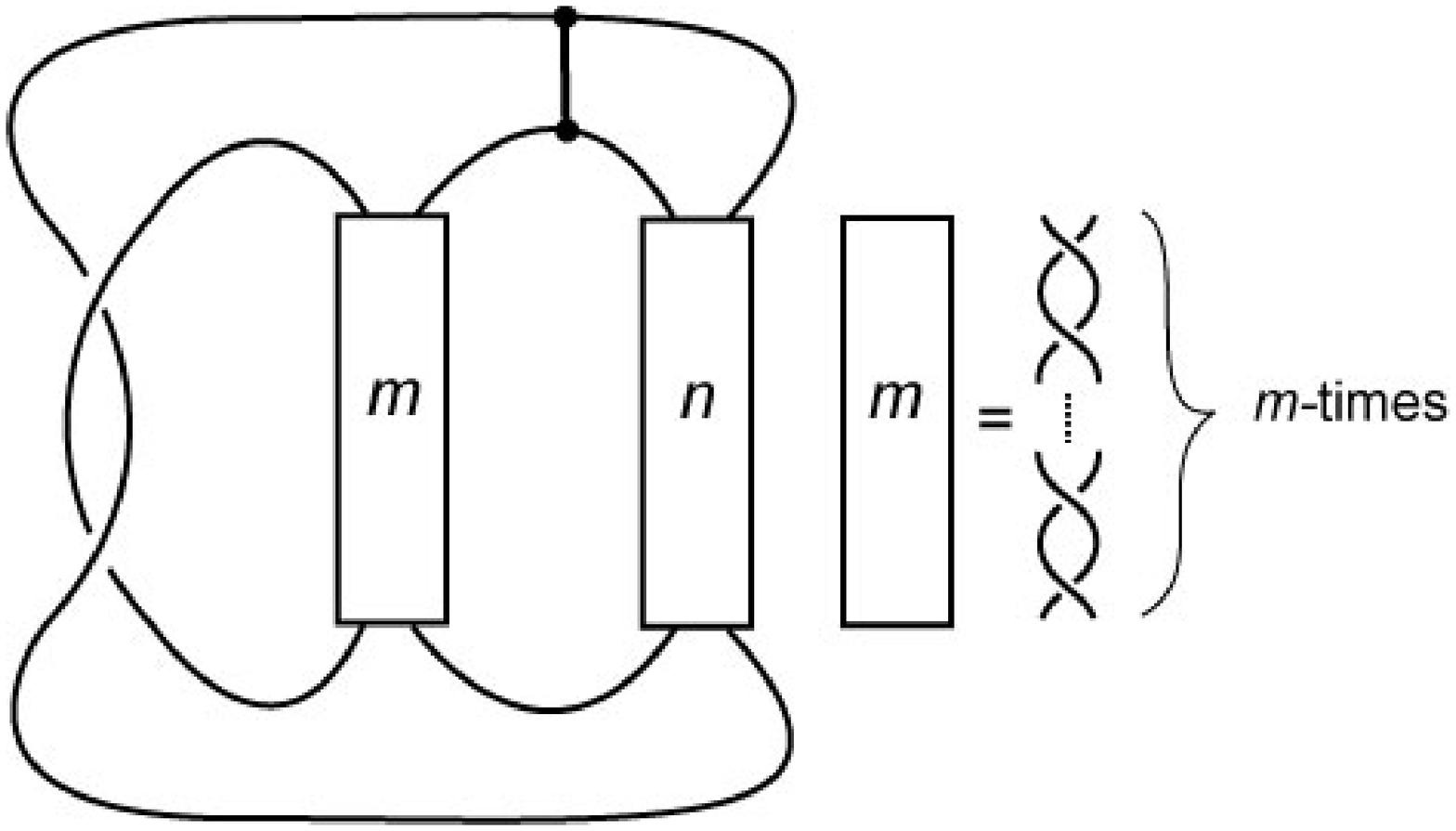}
\caption{}
\end{figure}

The pretzel knots of type $(-2,m,n)$ have projection
as illustrated in Figure 7.
All of them have a $(1,1)$-tunnel, for example,
we present a $(1,1)$-tunnel there.
Pretzel knots of this type 
have the $(1,1)$-tunnel in the same position.
Let $T$ be the unknotted torus as in Figure 8,
and we denote by $V_1$ (resp. $V_2$) the solid torus
$D^2\times S^1$ which is bounded `inside'
(resp. `outside') by $T$.
We can see that $(V_1,k_1)$ forms the pair of
a solid torus and a trivial arc.
We can also observe that $(V_2,k_2)$ is the
pair of them by an isotopy, which we call
$P$-isotopy.

\begin{figure}
\centering
\includegraphics[width=8.5cm,height=4.9cm]{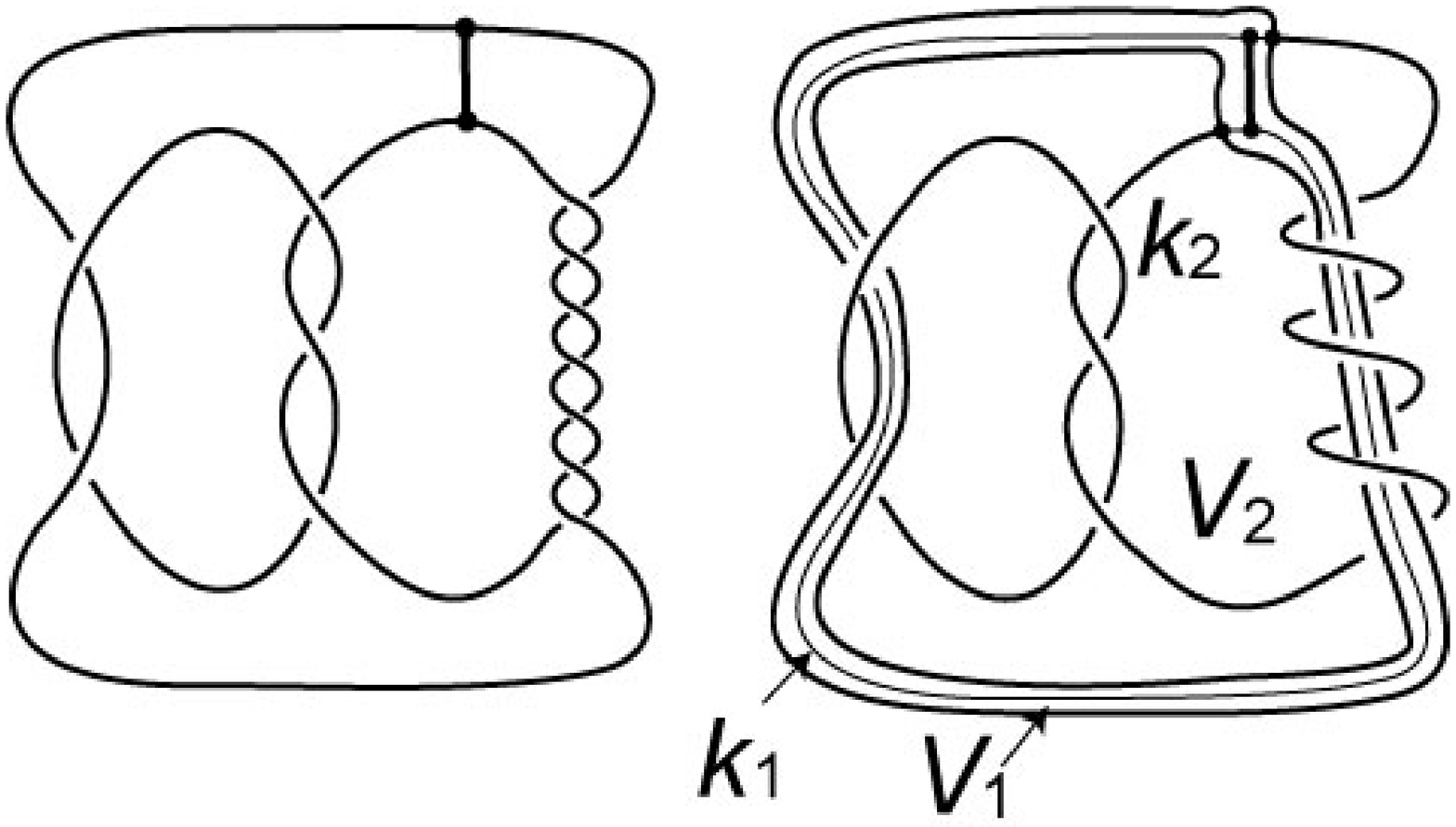}
\caption{}
\end{figure}

In what follows,
we consider the inverse of $P$-isotopy.

We can find a  meridian disk $D_2$ in
$V_2$ as in Figure 9.
We do proper isotopy for $k_2$ to
make $k$ an `almost' pretzel knot.
Here, 
since we need the information on 
$\bdd D_2$, we trace only $\bdd D_2$ under the
inverse of $P$-isotopy.

\begin{figure}
\centering
\includegraphics[width=8cm,height=3.4cm]{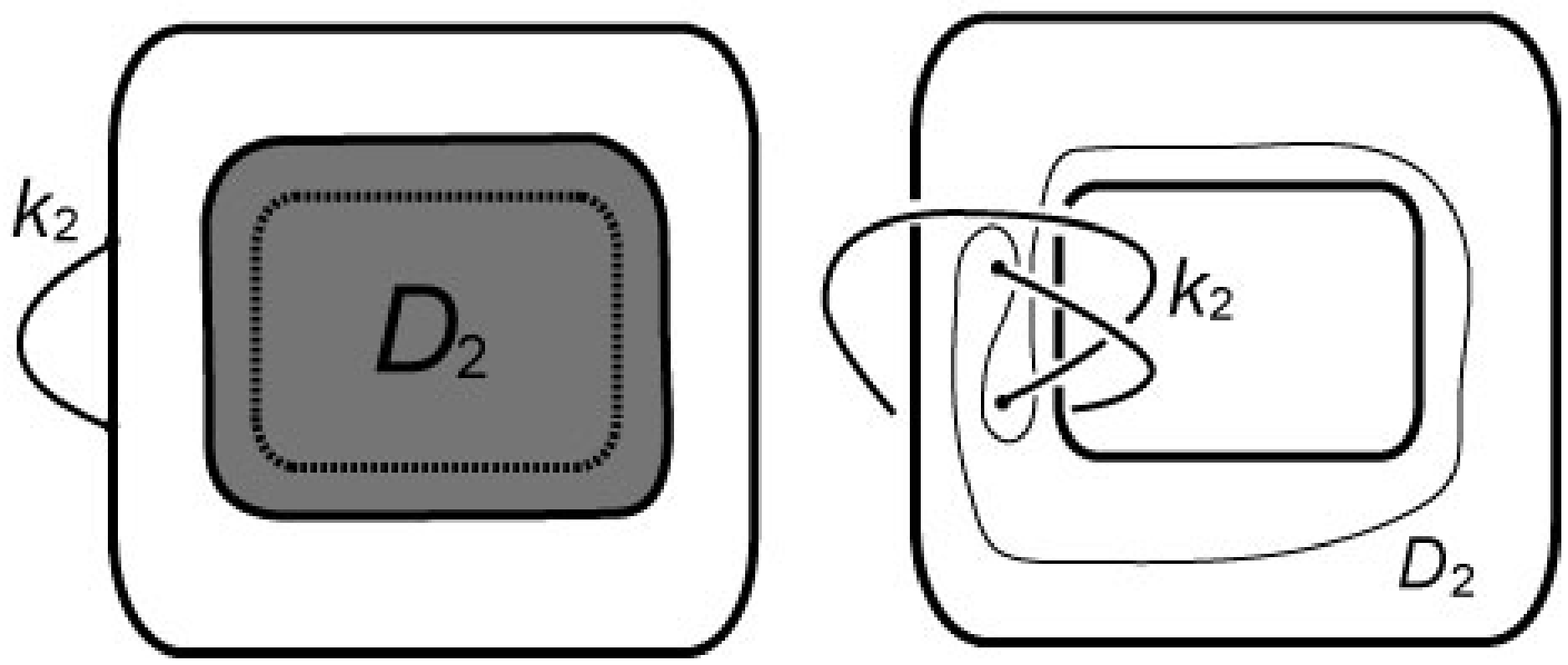}
\caption{}
\end{figure}

For the pretzel knot $P(-2,m,n)$,
we restore the part of `$-2$' and `$m$' part
as illustrated in Figures 9 and 10.
Continuing this way, we have the train track
as in the left-hand of Figure 11.

\begin{figure}
\centering
\includegraphics[width=8cm,height=3.4cm]{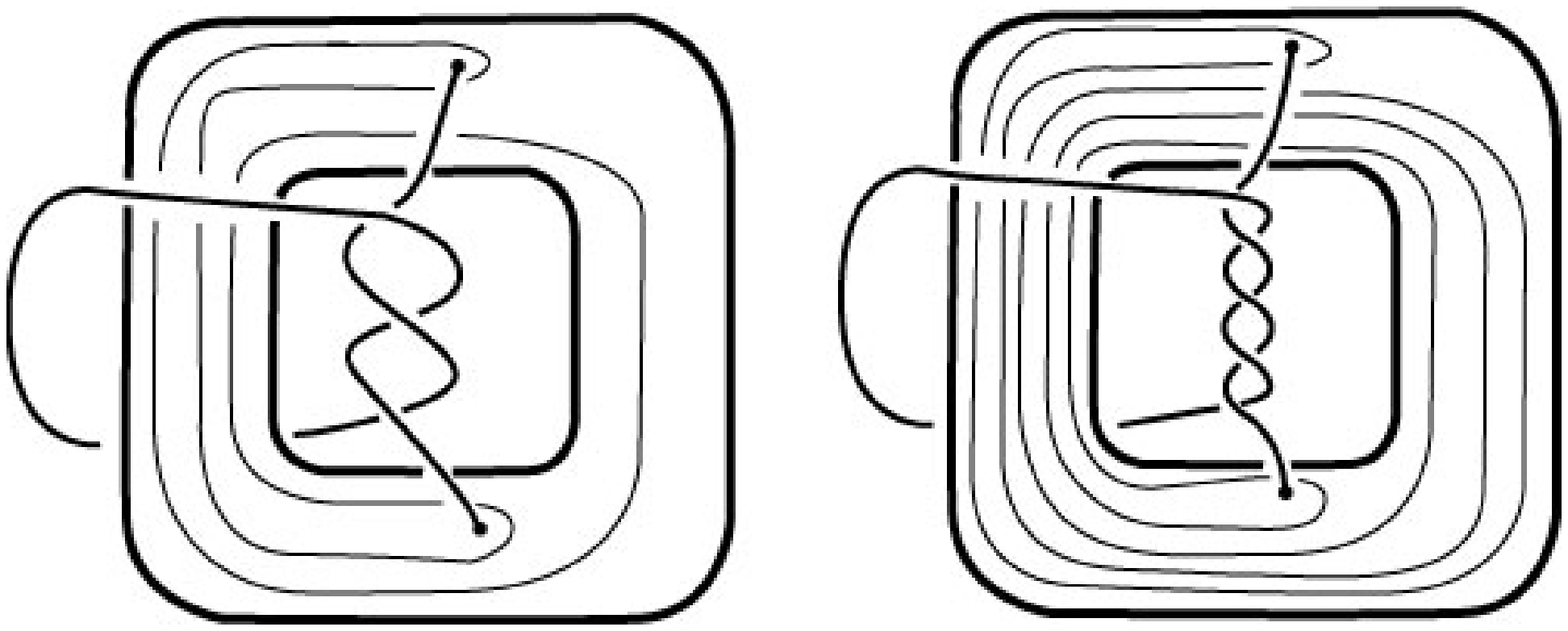}
\caption{}
\end{figure}

\begin{figure}
\centering
\includegraphics[width=8cm,height=3.4cm]{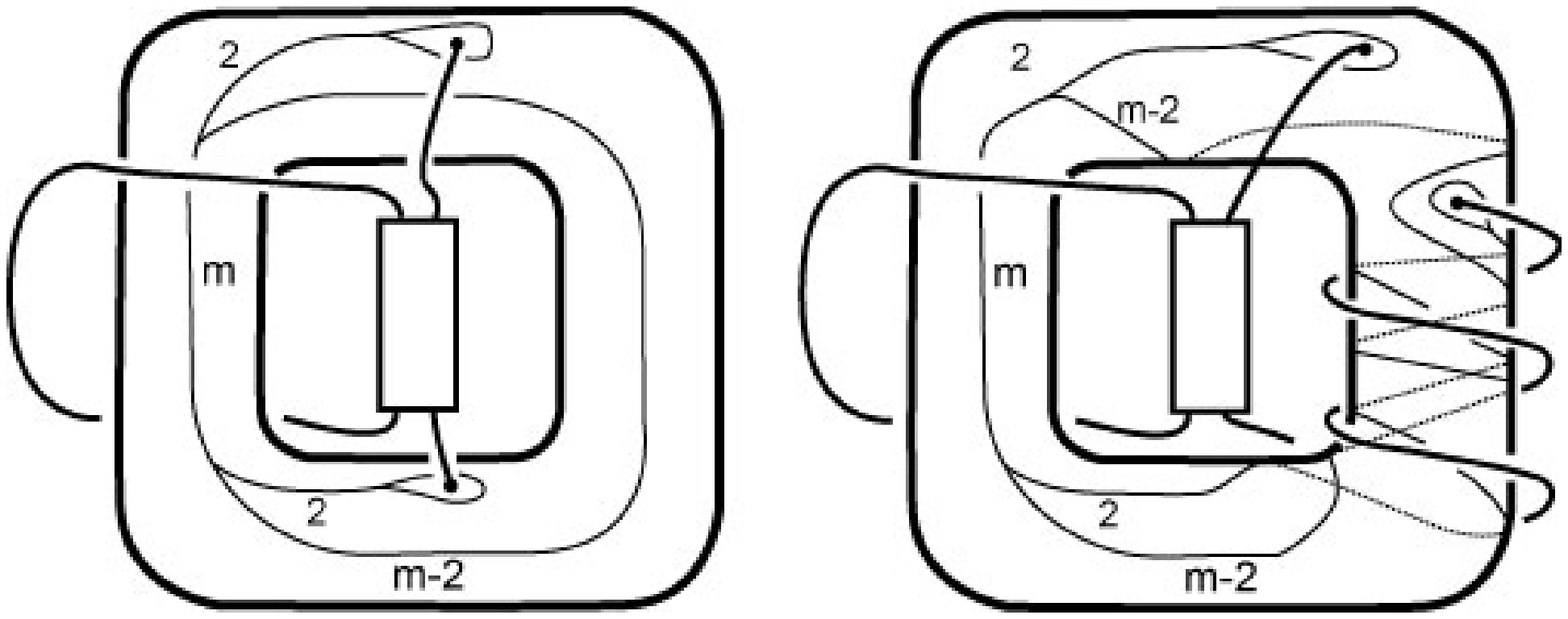}
\caption{}
\end{figure}

In order to restore the `$n$' part, 
we do the remaining 
of the inverse of $P$-isotopy 
as in right-hand of Figure 11.
Accordingly, 
we obtain a Heegaard diagram of the pretzel knot
of type $(-2,m,n)$. 
For example, 
a Heegaard diagram of the pretzel knot $P(-2,m,7)$ 
is illustrated
in Figure 12.

\begin{figure}
\centering
\includegraphics[width=9.7cm,height=5cm]{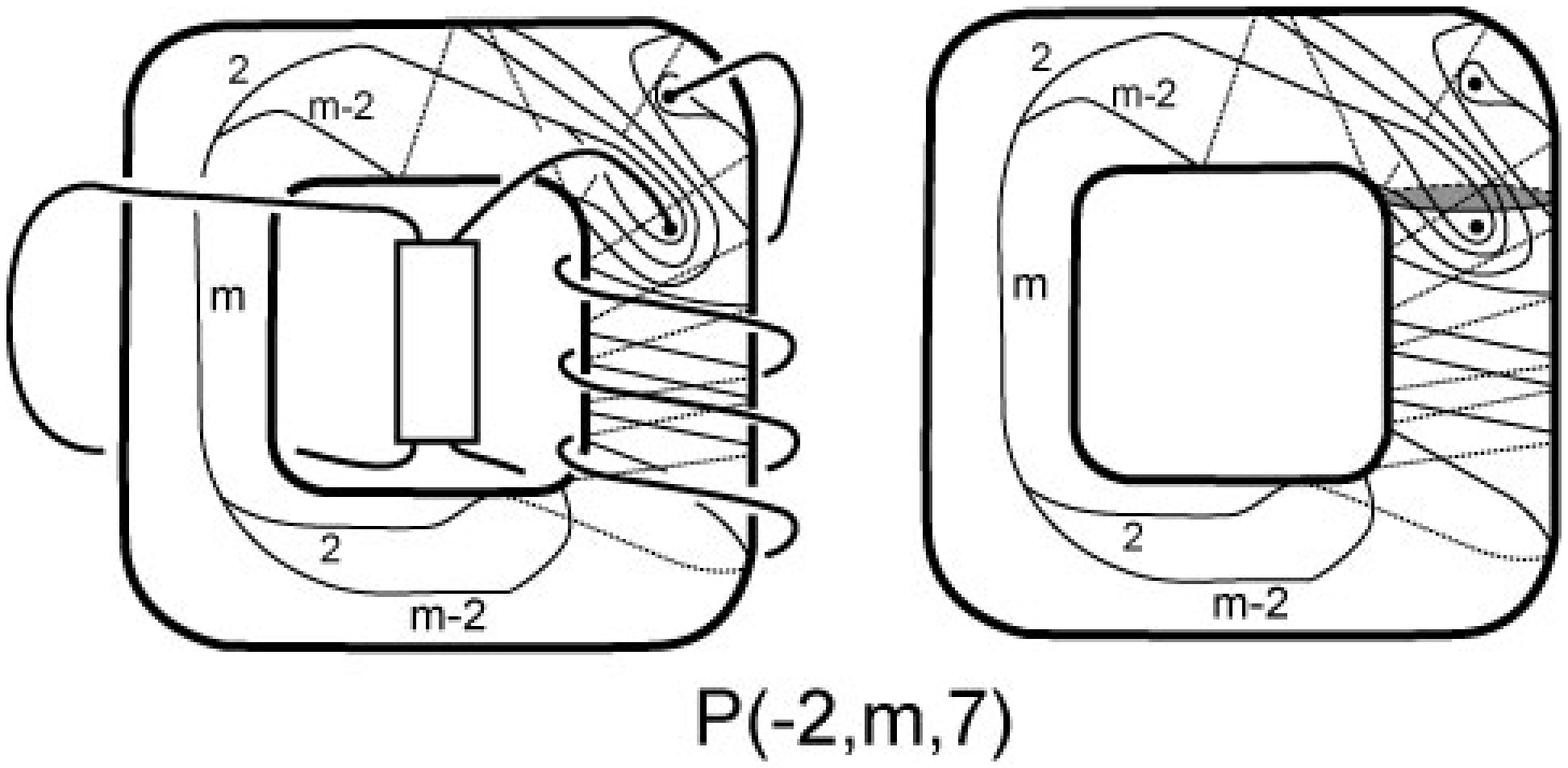}
\caption{}
\end{figure}

Therefore we have

\begin{lemma}
A Heegaard diagram for the pretzel knot $P(-2,m,n)$,
where $m$ and  $n$ are positive odd integers greater than one,
is carried by the train track
as illustrated in Figure $13$.
\end{lemma}

\begin{figure}
\centering
\includegraphics[width=8cm,height=8cm]{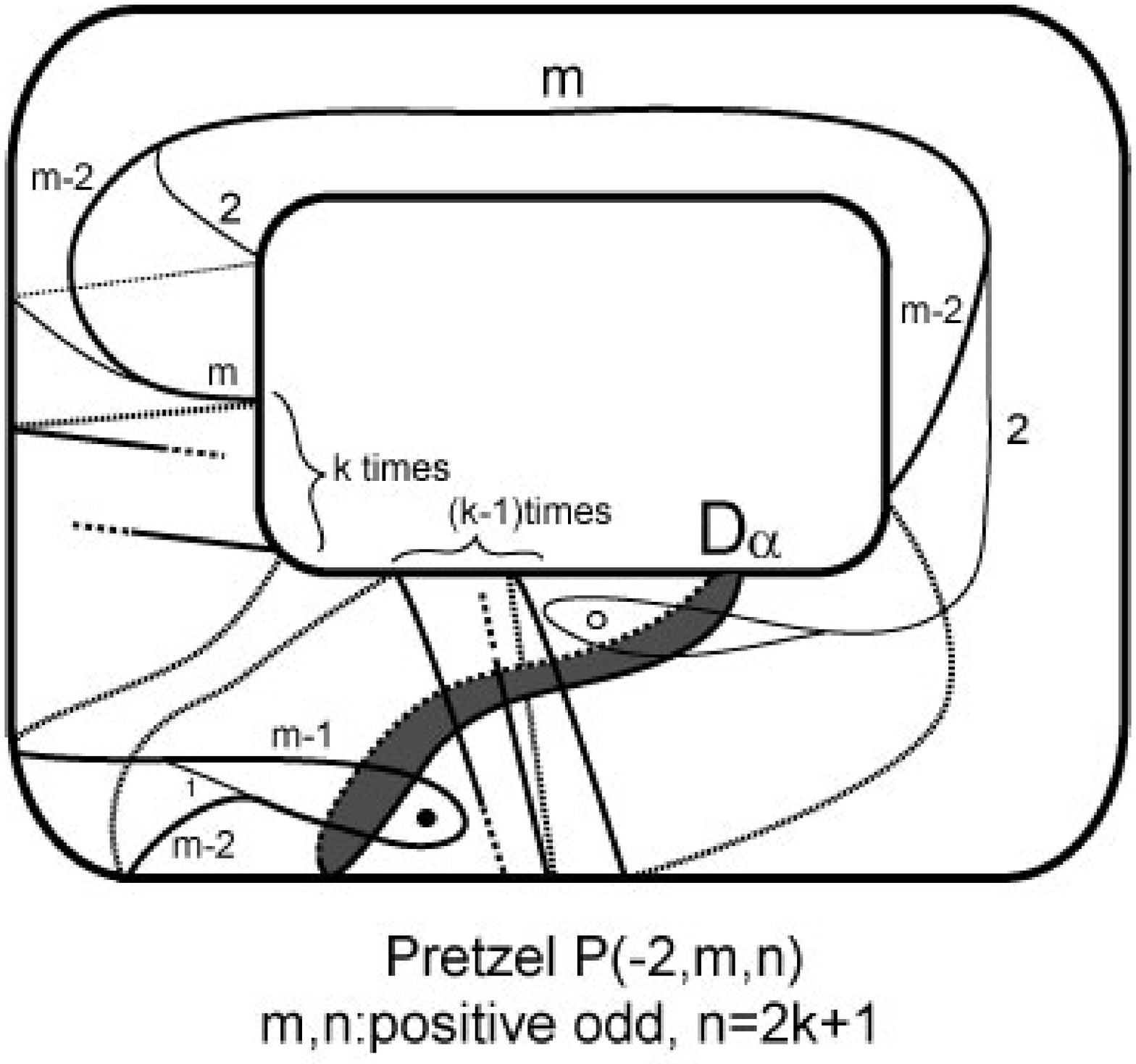}
\caption{}
\end{figure}


In order to treat uniformly, we assume 
$m\ge n\ge 5$ in the following.
After reading the proof for this, 
one can calculate the knot Floer homology
of the pretzel knots of type $(-2,m,3)$ 
by the same method.

By the straightforward calculation and taking a proper framing,
we have

\begin{lemma}
Suppose $n\ge 5$.
Figure $14$ illustrates the diagram
of $\widetilde{\alpha}$ and $\widetilde{\beta}$
on the universal cover $\C$ of the torus $\Sigma$
for the pretzel knot $P(-2,m,n)$,
where $\widetilde{\alpha}$ $($resp. $\widetilde{\beta})$
is a lift of $\partial D_\alpha^2$ $($resp. $\partial D_\beta^2)$
in $\C$.
\end{lemma}

\begin{figure}
\centering
\includegraphics[width=9.5cm,height=7.3cm]{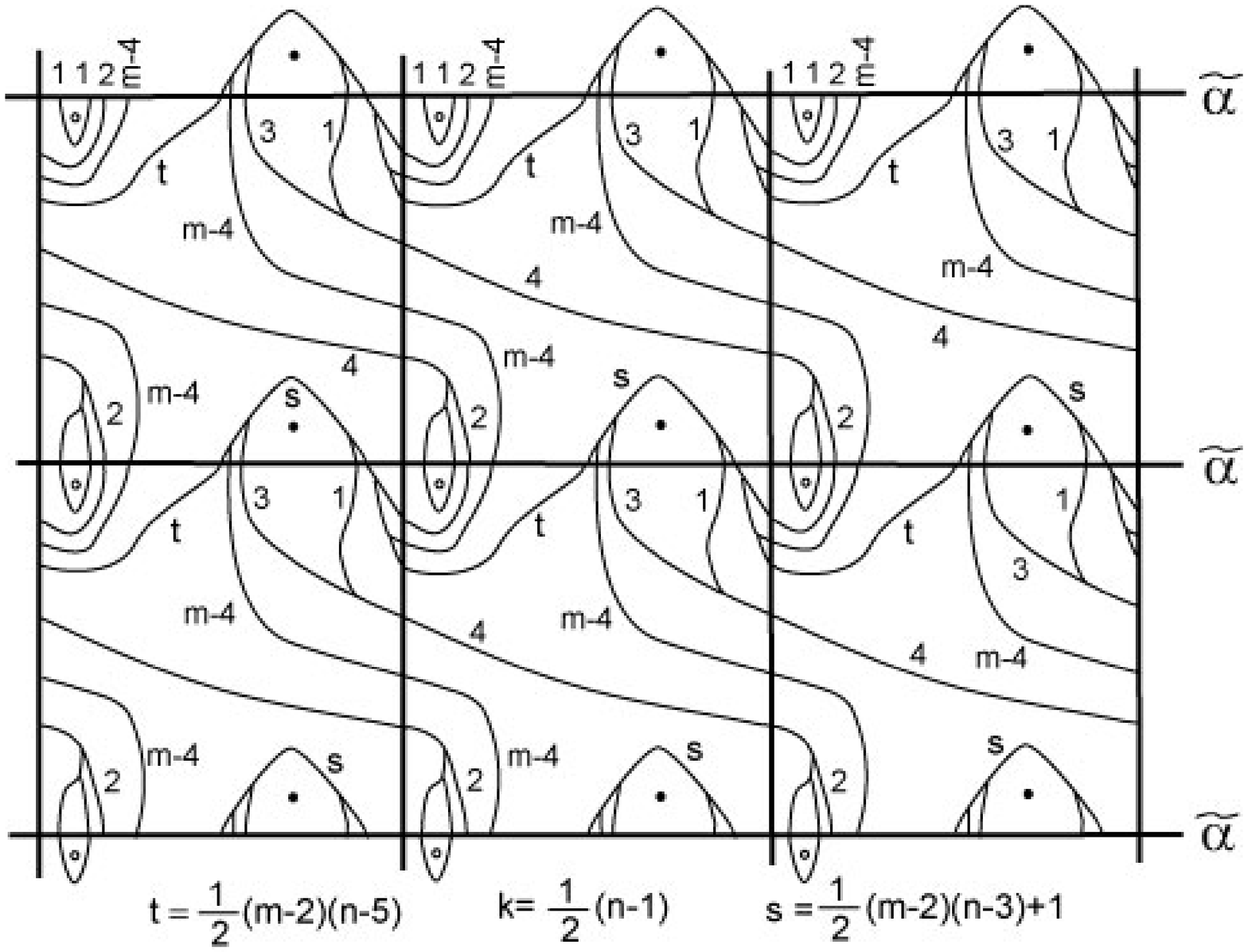}
\caption{}
\end{figure}

Next
we explicitly give 
the boundary operator from
the Heegaard diagram constructed above.
We use the following notations.
There is a holomorphic disc $\mathbb{D}_1$
(resp. $\mathbb{D}_2$) such that
it is bounded by an subarc in $\widetilde{\aaa}$
and that of $\widetilde{\bbb}$,
${\rm int}\mathbb{D}_1\cap(\widetilde{\aaa}\cup\widetilde{\bbb})=\emptyset$
(resp. ${\rm int}\mathbb{D}_2\cap(\widetilde{\aaa}\cup\widetilde{\bbb})=\emptyset$),
and $\mathbb{D}_1$ (resp. $\mathbb{D}_2$)
contains darkly-shaded (resp. lightly-shaded) point.
We denote $\widetilde{\aaa}\cap\widetilde{\bbb}\cap\bdd\mathbb{D}_1$
by $\widetilde{y}_1$ and $\widetilde{y}_2$ from the left hand side.
Similarly, let denote
$\widetilde{\aaa}\cap\widetilde{\bbb}\cap\bdd\mathbb{D}_2$
by $\widetilde{y}_3$ and $\widetilde{y}_4$.
We call $\widetilde{y}_1, \widetilde{y}_2,
\widetilde{y}_3, \widetilde{y}_4$
{\it exceptional generators\/} 
(see Figure 15).
Focus on $(m-2)$-lines of 
$\widetilde{\bbb}$ just the right hand side of 
$\T{y}_4$.
Toward left hand side, 
these lines run parallel 
until they intersect 
$(n-2)$ times with $\T{\alpha}$. 
We call these parallel lines a {\it bunch} of $(m-2)$-lines, 
and call $1$-{\it stage}, $2$-{\it stage},$\ldots$, 
$(n-2)$-{\it stage} from the left hand side at each part 
of the intersection of $\T{\alpha}$ 
and the bunch of $(m-2)$-lines. 
Thus 
each stage has $(m-2)$ points of 
$\T{\alpha}\cap\T{\beta}$, so we name 
$\T{x}_{i,1},\T{x}_{i,2},\ldots,\T{x}_{i,m-2}$ 
from the left hand side to the right hand side 
at the $i$-stage (see Figure 16).

\begin{figure}
\centering
\includegraphics[width=9.8cm,height=2.3cm]{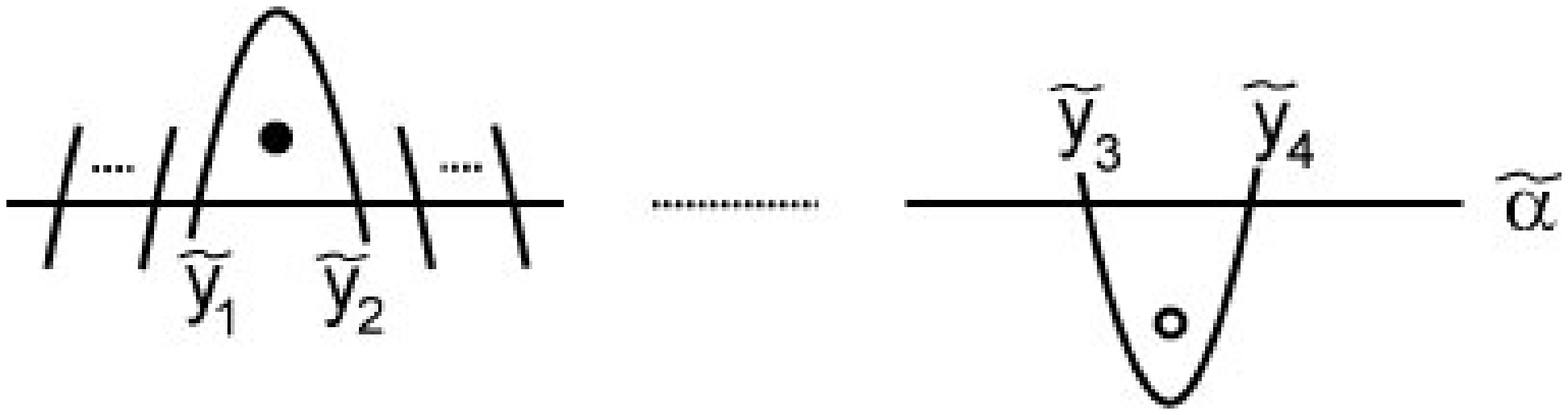}
\caption{}
\end{figure}

\begin{figure}
\centering
\includegraphics[width=9.5cm,height=6.3cm]{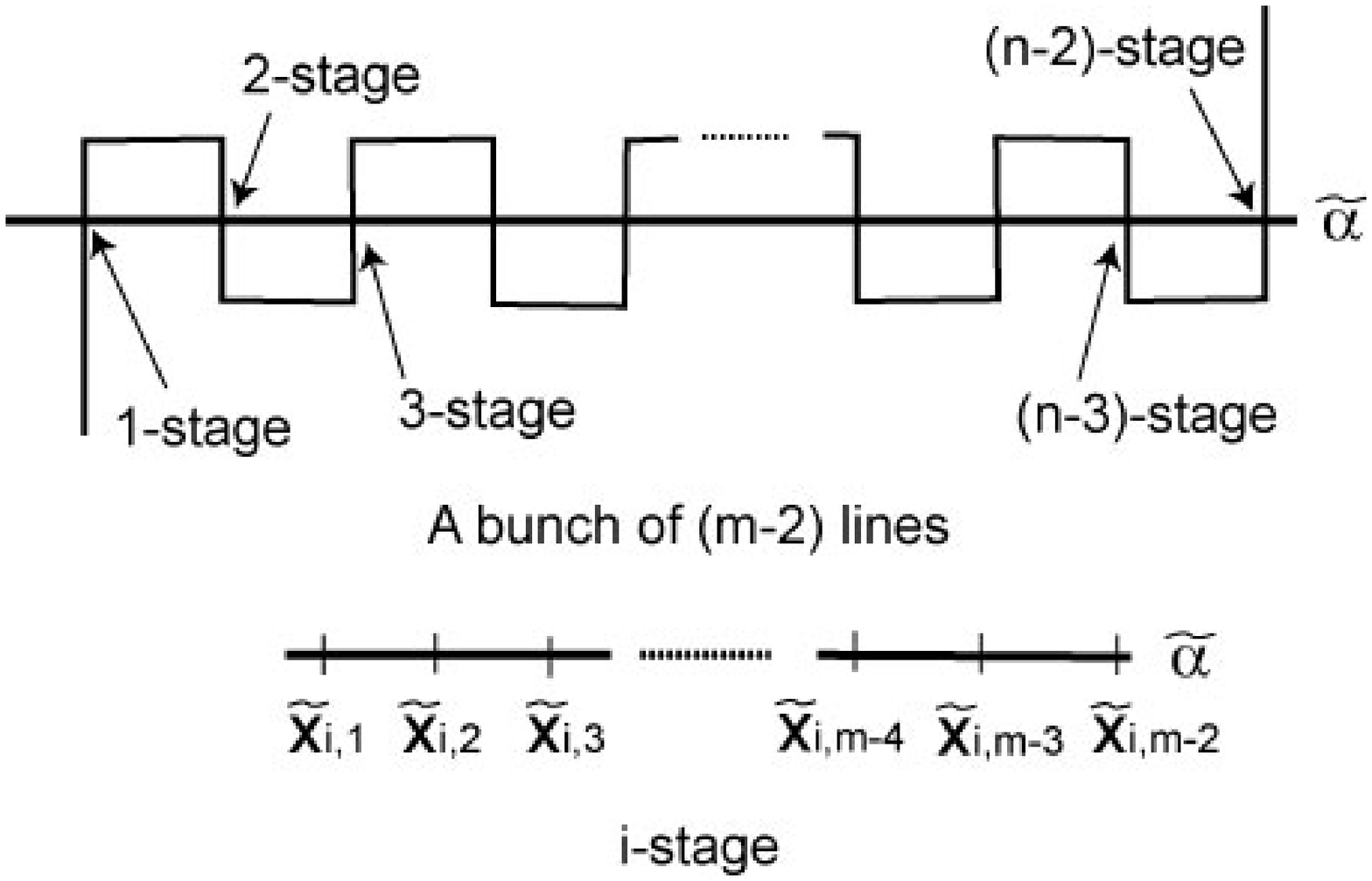}
\caption{}
\end{figure}

Then, 
by direct calculations as in Section 3, 
we see that 
$CFK^\infty\left(S^3,P(-2,m,n)\right)$ is 
generated as a $\Z$-module 
by generators indexed by $i\in\Z$ 
of the form:

\begin{align*}
&[\T{y}_1;i+\gamma(g)-1,i+\delta(g)+1],\quad 
[\T{y}_2;i+\gamma(g)-1,i+\delta(g)],\\
& [\T{y}_3;i+\delta(g),i+\gamma(g)-1],\quad 
[\T{y}_4;i+\delta(g)+1,i+\gamma(g)-1],\\
& [\T{x}_{2p+1,2q+1};i+\gamma(g)+p+q+1,i+\delta(g)-p-q]\quad
\left(0\leq p\leq n',\ 0\leq q\leq m'\right),\\
& [\T{x}_{2p+1,2q};i+\gamma(g)+p+q,i+\delta(g)-p-q]\quad
\left(0\leq p\leq n',\ 1\leq q\leq m'\right),\\
& [\T{x}_{2p,2q+1};i+\gamma(g)+m'+p-q,i+\delta(g)-m'-p+q]\quad
\left(1\leq p\leq n',\ 0\leq q\leq m'\right),\\
&[\T{x}_{2p,2q};i+\gamma(g)+m'+p-q,i+\delta(g)-m'-p+q-1]\quad
\left(1\leq p\leq n',\ 1\leq q\leq m'\right)
\end{align*}

\noindent
where 
$\displaystyle{
g=\frac{m+n}{2},\ m'=\frac{m-3}{2},\ n'=\frac{n-3}{2}
}$ and 
$$
\gamma(g)=
\left\{
\begin{gathered}
1-\frac{g-1}{2}\quad (g:\mathrm{odd})\\
1-\frac{g}{2}\ \qquad (g:\mathrm{even})
\end{gathered}
\right.
\qquad
\delta(g)=
\left\{
\begin{gathered}
\frac{g-1}{2}\qquad (g:\mathrm{odd})\\
\frac{g}{2}-1\qquad (g:\mathrm{even}).
\end{gathered}
\right.
$$
Then 
the boundary operator on 
$CFK^\infty\left(S^3,P(-2,m,n)\right)$ is 
explicitly given by the following lemmas.

\begin{lemma}
For exceptional generators: 
$\T{y}_1,\T{y}_2,\T{y}_3,\T{y}_4$, 
we have 

$\p [\T{y}_1;i,j]=[\T{y}_2;i,j-1],$

$\p [\T{y}_2;i,j]=0$,

$\p [\T{y}_3;i,j]=0$,

$\p [\T{y}_4;i,j]=-[\T{y}_3;i-1,j]$.
\end{lemma}

\begin{lemma}
For generators $\T{x}_{1,q}\ (1\leq q \leq m-2)$, we have 

$\p[\T{x}_{1,1};i,j]=[\T{x}_{2,m-2};i,j-1]+[\T{x}_{1,2};i,j-1]
-[\T{y}_2;i-2,j]$,

$\p[\T{x}_{1,2l};i,j]=[\T{x}_{2,m-2l-1};i,j-1]
\quad \left(1\leq l\leq m'\right)$,

$\p[\T{x}_{1,2l+1};i,j]=[\T{x}_{2,m-2l-2};i,j-1]
+[\T{x}_{1,2l+2};i,j-1]-[\T{x}_{1,2l};i-1,j]
\quad \left(1\leq l\leq m'-1\right)$,

$\p[\T{x}_{1,m-2};i,j]=[\T{x}_{2,1};i,j-1]-[\T{x}_{1,m-3};i-1,j]$.

\end{lemma}

\begin{lemma}
For generators $\T{x}_{n-2,q}\ (1\leq q \leq m-2)$, 
we have

$\p[\T{x}_{n-2,1};i,j]=[\T{x}_{n-2,2};i,j-1]-[\T{x}_{n-3,m-2};i-1,j]$,

$\p[\T{x}_{n-2,2l};i,j]=-[\T{x}_{n-3,m-2l-1};i-1,j]
\quad \left(1\leq l\leq m'\right)$,

$\p[\T{x}_{n-2,2l+1};i,j]=-[\T{x}_{n-3,m-2l-2};i-1,j]
-[\T{x}_{n-2,2l};i-1,j]$

$\qquad\qquad\qquad\qquad +[\T{x}_{n-2,2l+2};i,j-1]
\quad \left(1\leq l\leq m'-1\right)$,

$\p[\T{x}_{n-2,m-2};i,j]=
-[\T{x}_{n-3,1};i-1,j]-[\T{x}_{n-2,m-3};i-1,j]+[\T{y}_3;i,j-2]$.

\end{lemma}

\begin{lemma}
For generators $\T{x}_{2k,q}\ (1\leq k\leq n')$, 
we have

$\p[\T{x}_{2k,1};i,j]=[\T{x}_{2k,2};i-1,j]$,

$\p[\T{x}_{2k,2l};i,j]=0
\quad \left(1\leq l\leq m'\right)$,

$\p[\T{x}_{2k,2l+1};i,j]=[\T{x}_{2k,2l+2};i-1,j]-[\T{x}_{2k,2l};i,j-1]
\quad \left(1\leq l\leq m'-1\right)$,

$\p[\T{x}_{2k,m-2};i,j]=-[\T{x}_{2k,m-3};i,j-1]$.

\end{lemma}

\begin{lemma}
For generators $\T{x}_{2k+1,q}\ (1\leq k\leq n'-1)$, 
we have

$\p[\T{x}_{2k+1,1};i,j]=[\T{x}_{2k+2,m-2};i,j-1]+[\T{x}_{2k+1,2};i,j-1]
-[\T{x}_{2k,m-2};i-1,j]$,

$\p[\T{x}_{2k+1,2l};i,j]=[\T{x}_{2k+2,m-2l-1};i,j-1]-[\T{x}_{2k,m-2l-1};i-1,j]
\quad \left(1\leq l\leq m'\right)$,

$\p[\T{x}_{2k+1,2l+1};i,j]
=[\T{x}_{2k+2,m-2l-2};i,j-1]+[\T{x}_{2k+1,2l+2};i,j-1]$

$\qquad\qquad\qquad\quad -[\T{x}_{2k,m-2l-2};i-1,j]-[\T{x}_{2k+1,2l};i-1,j]
\quad \left(1\leq l\leq m'-1\right)$,

$\p[\T{x}_{2k+1,m-2};i,j]
=[\T{x}_{2k+2,1};i,j-1]
-[\T{x}_{2k+1,m-3};i-1,j]-[\T{x}_{2k,1};i-1,j]$.

\end{lemma}

By using these lemmas, 
we can immediately obtain Theorem 5.1. 
In particular, 
we easily see that 
the Floer homology class of 
$\widehat{HF}(S^3)\cong\Z$ is 
represented by the generator $[\T{y}_4;0,-g]$ 
(so its absolute grading is zero). 
Accordingly, 
we can also conclude that 
the invariant $\tau(P(-2,m,n))$ is 
equal to $-g$.

\bigskip 

\begin{figure}
\centering
\includegraphics[width=9.5cm,height=19.5cm]{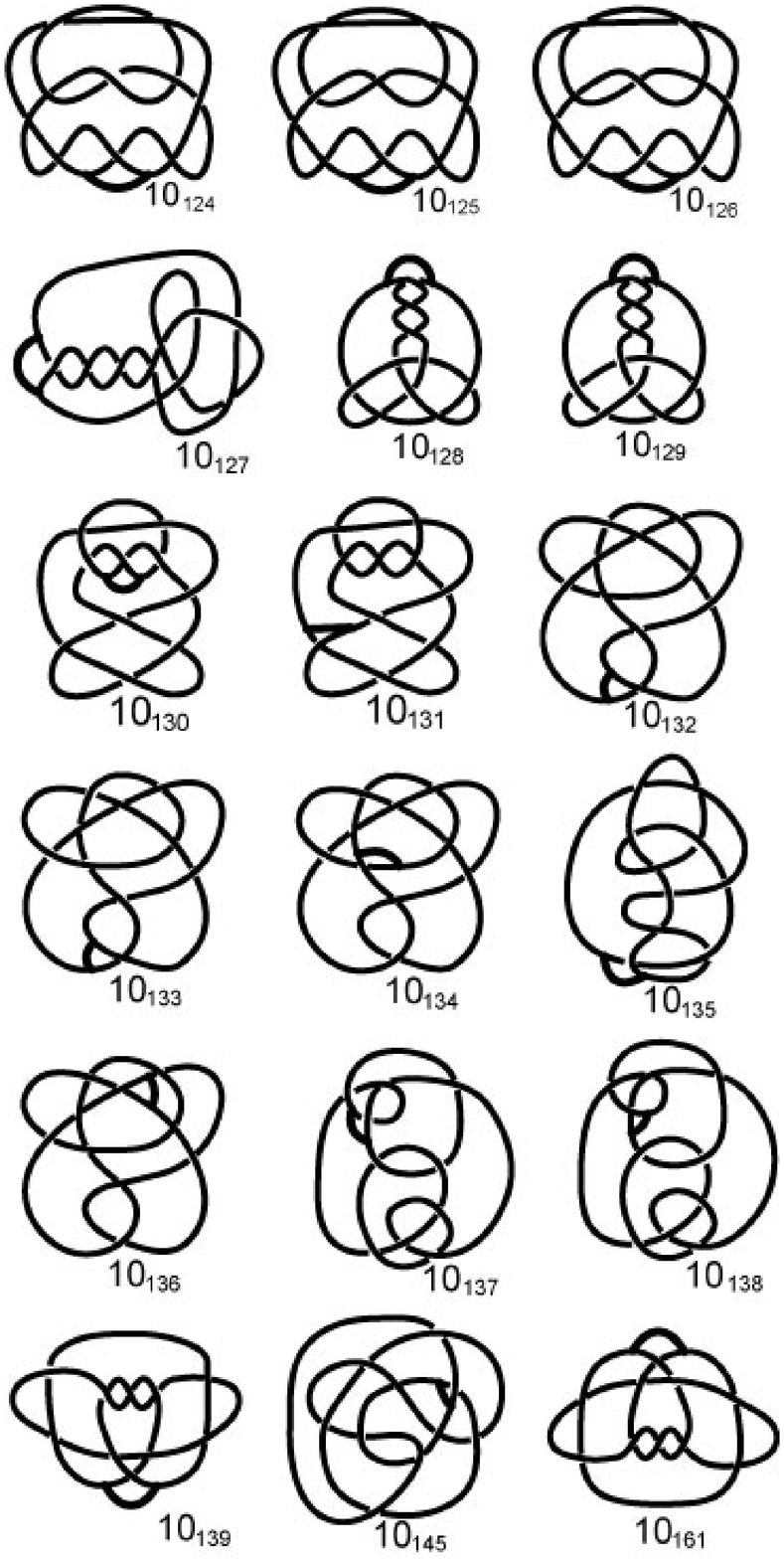}
\caption{}
\end{figure}








\end{document}